\documentclass[a4paper,oneside,12pt]{article}
\usepackage{amssymb,amsmath,amsfonts,amsthm}
\usepackage{amsmath}
\usepackage{graphicx,graphics}
\usepackage{epsfig}
\usepackage{latexsym}
\usepackage[applemac]{inputenc}
\usepackage{ae,aecompl}
\usepackage{mathrsfs}
\usepackage[english]{babel}
 \usepackage[colorlinks=true]{hyperref}
\usepackage{pstricks}

\renewcommand{\leq}{\lesqlant}

\newcommand{\E}[1]{\mathbb{E} \left [ #1 \right ]}

\newcommand{\R}{\mathbb R}

\newcommand{\Z}{\mathbb Z}

\newcommand{\ligne}{\vspace{3mm}}

\newsavebox{\fmbox}

\renewcommand{\geq}{\geqslant}
\renewcommand{\leq}{\leqslant}

\newcommand{\op}[1]{\operatorname{#1 }}

\newcommand{\tends}[2]{\overset{#1}{\underset{#2}{\longrightarrow}}}

\renewcommand{\tends}[2]{\overset{#1}{\underset{#2}{\longrightarrow}}}

\newtheorem{theorem}{Theorem}[section]

\newtheorem{proposition}[theorem]{Proposition}
\newtheorem{lemma}[theorem]{Lemma}
\newtheorem{corollary}[theorem]{Corollary}
\newtheorem{question}{Question}
\newtheorem{definition}[theorem]{Definition}
\newtheorem{exe}{Example}

\newtheorem{rek}[theorem]{Remark}

\usepackage{enumerate}

\newcommand{\Eb}[1]{\mathbf{E}\left[#1\right]}

\date{}
\title{Ergodic Theory on Stationary Random Graphs}
\author{Itai Benjamini and  Nicolas Curien}
\begin{document}
\maketitle

\begin{abstract}A stationary random graph is a random rooted graph whose distribution is invariant under re-rooting along the simple random walk. We adapt the entropy technique developed for Cayley graphs and show in particular that stationary random graphs of  subexponential growth are almost surely Liouville, that is, admit no non constant bounded harmonic functions. Applications include the uniform infinite planar quadrangulation and  long-range percolation clusters. \end{abstract}

\section{Introduction}

A stationary random graph $(G,\rho)$ is a random rooted  graph whose distribution is invariant under re-rooting along a simple random walk started at the root $\rho$ (see Section \ref{defi} for a precise definition). The entropy technique and characterization of the Liouville property for groups,  homogeneous graphs or random walk in random environment  \cite{Kai90,Kai98,Kai01,KKR04,KV83,KW02} are adapted  to this context. In particular we have
\begin{theorem}\label{subexp} Let $(G,\rho)$ be a stationary random graph of \emph{subexponential growth} in the sense that  \begin{eqnarray} n^{-1}\mathbf{E}\Big[\log\big(\#B_{G}(\rho,n)\big)\Big] &\underset{n \to \infty}{\longrightarrow}&0 \label{subexpo} ,\end{eqnarray} where $\#B_G(\rho,n)$ is the number of \emph{vertices} within distance $n$ from the root $\rho$, then $(G,\rho)$ is almost surely Liouville.
\end{theorem}
 Recall that a function from the vertices of a graph to  $\mathbb{R}$ is harmonic if and only if the value of the function at a vertex is the average of the value over its neighbors, for all vertices of the graph. We call graphs admitting no non constant bounded  harmonic functions {\em Liouville}. In the case of graphs of bounded degree we show in Proposition \ref{fundamental} that stationary non-Liouville random graphs are \emph{ballistic}.  \medskip

One of the motivation of this work lies in the study of the Uniform Infinite Planar Quadrangulation (abbreviated by UIPQ)  introduced  in \cite{Kri05} (following the pioneer work of \cite{AS03}). The UIPQ  is a stationary random infinite planar graph whose faces are all squares. This object  is very natural and of special interest for understanding two dimensional quantum gravity and has triggered a lot of work, see e.g. \cite{Ang03,AS03,CD06,CMMinfini,LGM10,Men08}. One of the  fundamental questions regarding the UIPQ, is to prove recurrence or transience of simple random walk on this graph. Unfortunately, the degrees in the UIPQ are not bounded thus the techniques of \cite{BS01} fail to apply. Nevertheless it has been conjectured in \cite{AS03} that the UIPQ is a.s.\,recurrent. As an application of Theorem \ref{subexp}, we deduce a step in this direction, 
\begin{corollary} \label{uipq} The Uniform Infinite Planar Quadrangulation is almost surely Liouville.
\end{corollary}
See also the very recent work of Steffen Rohde and James T.\,Gill \cite{GR10} proving that the conformal type of the Riemann surface associated to the UIPQ is parabolic. Another application concerns a question of Berger  \cite{Berg02} and consists in proving that certain long range percolation clusters are Liouville (see Section 5.2).
\medskip

The notion of stationary random graph generalizes the concepts of Cayley and transitive graph where the homogeneity of the graph is replaced by stationarity along the simple random walk. This notion is very closely related to the ergodic theory notions of unimodular random graphs of \cite{AL07} and measured equivalence relations see e.g.\,\cite{Kai98,Kai01,KS10,Pau99}. Roughly speaking, unimodular random graphs correspond, after biasing by the degree of the root, to stationary and \emph{reversible} random graphs (see Definition \ref{seuldef}). We then reinterpret  ideas from measured equivalence relations theory to prove (Theorem \ref{ballistic}) that if a stationary random graph of bounded degree $(G,\rho)$ is non reversible then the simple random walk on $G$ is ballistic, thus improving Theorem A of \cite{Pau99} and extending \cite{SW90} in the case of transitive graphs.\\

In \cite{BCsnake} the authors also use the notions of stationary and unimodular random graph  in order to show that the simple random walk on $\mathbb{Z}^d$ indexed by $T_{\infty}$, the critical geometric Galton-Watson tree conditioned to survive \cite{Kes86}, is recurrent if and only if $d \leq 4$.\\

The goal of this paper is not to prove striking new results, indeed much of the general results stated in this work are adaptations or variants of known results in the context of measured equivalence relations. Rather, we present them in a new and clear probabilistic framework which is of independent interest.  We thus chose to focus on the Liouville property for graphs and its application to the UIPQ as main direction. However it is believable that larger parts of the theory of equivalence relations can be adapted to the random graph setting.  In the last section, we also construct (Proposition \ref{counterexample})  a stationary and reversible random graph of subexponential growth which is planar and transient. This indicates that the theory of local limits of random planar graphs of bounded degree developped in \cite{BS01} can not be extended to the unbounded degree case in a straightforward manner. \\

The paper is organized as follows. The remainder of this section is devoted to a formal definition of stationary and reversible random graphs. Section 2 recalls the links between these concepts, unimodular random graphs and measured equivalence relations. The entropy technique is developed in Section 3. In Section 4 we explore under which conditions a stationary random graph is not reversible. The last section is devoted to applications and open problems. \medskip

\subsection{Definitions}
\label{defi}
A graph $G=(\op{V}(G),\op{E}(G))$ is a pair of sets, $\op{V}(G)$ representing the set of vertices and $\op{E}(G)$ the set of (unoriented) edges. In the following,  all the graphs considered are countable, connected and locally finite. We also restrict ourself to \emph{simple} graphs, that is, without loops nor multiple edges. 
Two vertices $x,y \in \op{V}(G)$ linked by an edge are called \emph{neighbors} in $G$ and we write $x\sim y$. The \emph{degree} $\op{deg}(x)$ of $x$ is the number of neighbors of $x$ in $G$.  For any pair $x,y \in G$, the \emph{graph distance} $\op{d}_{\op{gr}}^{G}(x,y)$ is the minimal length of a path joining $x$ and $y$ in $G$. For every $r \in \mathbb{Z}_{+}$, the ball of radius $r$ around $x$ in $G$ is the subgraph of $G$ spanned by the vertices at distance less than or equal to $r$ from $x$ in $G$, it is denoted by $B_{G}(x,r)$. \medskip

A \emph{rooted graph} is a pair $(G,\rho)$ where $\rho\in \op{V}(G)$ is called the root vertex. An isomorphism between two rooted graphs is a graph isomorphism that maps the roots of the graphs. Let $\mathcal{G}_{\bullet}$ be the set of isomorphism classes of locally finite rooted graphs $(G,\rho)$, endowed with the distance $\op{d_{loc}}$ defined by
$$ \op{d_{loc}}\big((G_{1},\rho_{1}),(G_{2},\rho_{2})\big) = \inf\left \{ \frac{1}{r+1} : r \geq 0 \mbox{ and } (B_{G_{1}}(\rho_{1},r),\rho_{1}) \simeq (B_{G_{2}}(\rho_{2},r),\rho_{2})\right\},$$ where $\simeq$ stands for the rooted graph equivalence. With this topology, $\mathcal{G}_{\bullet}$ is a Polish space (see \cite{BS01}). Similarly, we define  $\mathcal{G}_{\bullet \bullet}$ (resp.\,$\vec{\mathcal{G}}$) to be the set of isomorphism classes of bi-rooted graphs $(G,x,y)$ that are graphs with two distinguished ordered points (resp.\,graphs $(G,(x_{n})_{n\geq0})$ with a semi-infinite path), where the isomorphisms considered have to map the two distinguished points (resp.\,the path).  These two sets are equipped with variants of the distance $\op{d_{loc}}$ and are Polish with the induced topologies.   Formally elements of $\mathcal{G}_{\bullet}, \mathcal{G}_{\bullet\bullet}$ and $\vec{\mathcal{G}}$ are equivalence classes of graphs, but we will not distinguish between graphs and their equivalence classes and we use the same terminology and notation. One way to bypass this identification is to choose once for all a canonical representative in each class, see \cite[Section 2]{AL07}.

Let $(G,\rho)$ be a rooted graph.  For $x \in \op{V}(G)$ we denote the law of the \emph{simple random walk} $(X_{n})_{n\geq 0}$ on $G$ starting from $x$ by $\mathrm{P}_{x}^G$ and its expectation by $\mathrm{E}_{x}^G$. Formally this makes no sense since $(G,\rho)$ is an equivalence class of graphs, however it is easy to check that the distribution of $(G,(X_{n})_{n\geq 0}) \in \vec{\mathcal{G}}$ when $(X_{n})$ starts from $\rho$ is well-defined, that is does not depend on the representative chosen for $(G,\rho)$. We speak of ``the simple random walk of law $\mathrm{P}_{\rho}^G$ conditionally on $(G,\rho)$''. It is easy to check that all the quantities we will use in the paper do not depend of a choice of a representative of $(G,\rho)$.

A random rooted graph $(G,\rho)$ is a random variable taking values in $\mathcal{G}_{\bullet}$. In this work we will use $\mathbf{P}$ and $\mathbf{E}$ for the probability and expectation referring to the underlying random graph. If conditionally on $(G,\rho)$, $(X_{n})_{n\geq 0}$ is the simple random walk started at $\rho$, we denote the distribution of $(G, (X_{n})_{n\geq 0}) \in \vec{\mathcal{G}}$ by
$ \mathbb{P}$, and  the corresponding expectation by $\mathbb{E}$. The following concept is quite standard. 
\begin{definition} \label{seuldef}Let $(G,\rho)$ be a random rooted graph. Conditionally on $(G,\rho)$, let $(X_{n})_{n\geq 0}$ be the simple random walk on $G$ starting from $\rho$. The graph $(G,\rho)$ is called \emph{stationary} if
\begin{eqnarray}
(G,\rho) &= &(G,X_{n}) \label{sta} \quad \mbox{in distribution, for all }n\geq 1,
\end{eqnarray}
or equivalently for $n=1$. In words a stationary random graph is a random rooted graph whose distribution is invariant under re-rooting along a simple random walk on $G$. Furthermore, $(G,\rho)$ is called  \emph{reversible} if
\begin{eqnarray}
(G,X_{0},X_{1}) &= &(G,X_{1},X_{0}) \label{reversible}\quad \mbox{in distribution.}
\end{eqnarray}

\end{definition}

Clearly any reversible random graph is stationary. Note that our reversibility condition is different from the usual notion for Markov processes.

\begin{exe} Any Cayley graph rooted at any vertex is stationary and reversible. Any transitive graph $G$ (i.e. whose isomorphism group is transitive on $\op{V}(G)$) is stationary.  For examples of transitive graphs which are not reversible, see \cite[Examples 3.1 and 3.2]{BLPS99}. E.g.\,the ``grandfather'' graph (see Fig.\,below) is a transitive (hence stationary) graph which is not reversible. \begin{center} \includegraphics[height=6cm]{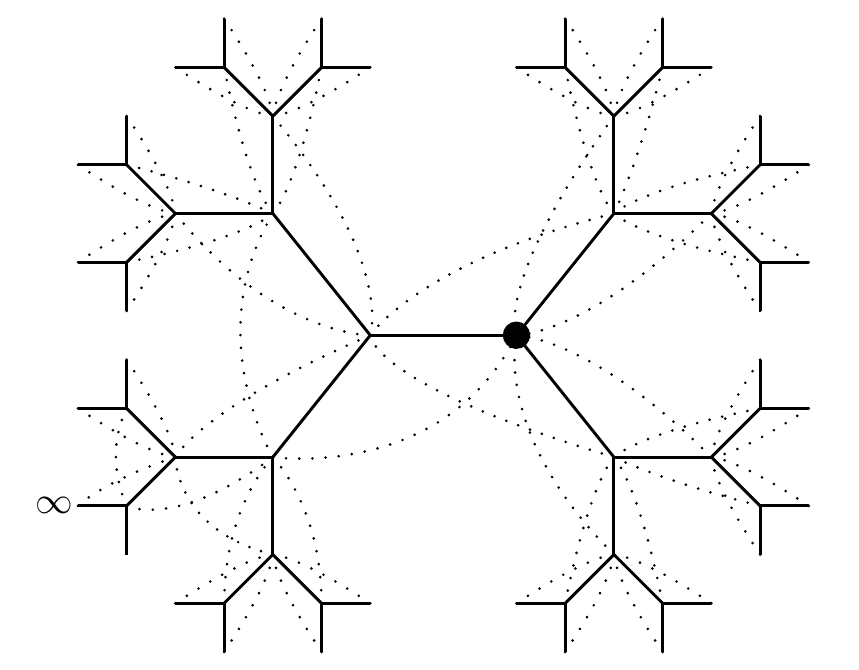}\\
 \textsc{Fig.: }The ``grandfather'' graph is obtained from the 3-regular tree by choosing a point at infinity that orientates the graph and adding all the edges from grandsons to grandfathers. \end{center}
\end{exe}
\begin{exe} \label{example} \cite[Section 3.2]{BS01} Let $G$ be a finite connected graph. Pick a vertex $\rho \in \op{V}(G)$ with a probability proportional to its degree (normalized by $\sum_{u \in \op{V}(G)}\op{deg}(u)$). Then $(G,\rho)$ is a reversible random graph.
\end{exe}
\begin{exe}[Augmented Galton-Watson tree] \label{AGW} Consider two independent Galton-Watson trees with offspring distribution $(p_{k})_{k \geq 0}$. Link the root vertices of the two trees by an edge and root the obtained graph at the root of the first tree. The resulting random rooted graph  is stationary and reversible, see \cite{KS10,LPP95,LP10}.
\end{exe}
\section{Connections with other notions}
As we will see, the concept of stationary random graph can be linked to various notions. In the context of bounded degree, stationary random graphs generalize unimodular random graphs \cite{AL07}. Stationary random graphs are closely related to graphed equivalence relations with a harmonic measure, see \cite{Kai01,Pau99}. We however think that the probabilistic Definition \ref{seuldef} is more natural for our applications.
\label{links}

\subsection{Ergodic theory} \label{ergo}
We formulate  the notion of  stationary random graphs in terms of ergodic theory.
We can define the \emph{shift} operator $\theta$ on $ \vec{\mathcal{G}}$ by $ \theta\big((G,(x_{n})_{n\geq 0})\big) = \big(G,(x_{n+1})_{n\geq 0}\big),$ and the projection $\pi: \vec{\mathcal{G}}\to\mathcal{G}_{\bullet}$ by $\pi\big((G,(x_{n})_{n\geq0})\big)= (G,x_{0}).$

Recall from the last section that if $\mathbf{P}$ is the law of $(G,\rho)$ we write $\mathbb{P}$ for the distribution of $(G,(X_{n})_{n\geq0})$ where $(X_{n})_{n\geq 0}$ is the simple random walk on $G$ starting at $\rho$.
The following proposition is a straightforward translation of the notion of a stationary random graph into that of a $\theta$-invariant probability measure on $\vec{\mathcal{G}}$.
\begin{proposition}\label{correspondence} Let $\mathbf{P}$ a probability measure on $\mathcal{G}_{\bullet}$ and $\mathbb{P}$ the associated probability measure on $\vec{\mathcal{G}}$. Then $\mathbf{P}$ is stationary if and only if $\mathbb{P}$ is invariant under $\theta$.
\end{proposition}
As usual, we will say that $\mathbb{P}$ (and by extension $\mathbf{P}$ or directly $(G,\rho)$) is ergodic if $\mathbb{P}$  is ergodic for $\theta$. Proposition \ref{correspondence} enables us to use all the powerful machinery of ergodic theory in the context of stationary random graphs. For instance, the classical theorems on the range and speed of a random walk on a group are valid:

\begin{theorem} \label{range} Let $(G,\rho)$ be a stationary and ergodic random graph. Conditionally on $(G,\rho)$ denote $(X_{n})_{n\geq0}$ the simple random walk on $G$ starting from $\rho$. Set $R_{n} = \# \{X_{0}, \ldots , X_{n}\}$ and $D_{n}=\op{d}_{\op{gr}}^G(X_{0}, X_{n})$ for the range and distance from the root of the random walk at time $n$. There exists a constant $s\geq 0$ such that we have the following almost sure and $\mathbb{L}^1$ convergences for $\mathbb{P}$,
\begin{eqnarray} \frac{R_{n}}{n} &\tends{a.s. \ \mathbb{L}^1}{n \to \infty}& \mathbb{P}\left( \bigcap_{i \geq 1} \{ X_{i} \ne \rho \} \right), \label{range1}\\
\frac{D_{n}}{n} &\tends{a.s. \ \mathbb{L}^1}{n\to \infty}& s.\label{speed} \end{eqnarray}\end{theorem}
\begin{rek} \label{transience} In particular a stationary and ergodic random graph is transient if and only if the range of the simple random walk on it grows linearly.
\end{rek}
\proof The two statements are straightforward adaptations of \cite{Der80}. See also \cite[Proposition 4.8]{AL07}. \endproof


\subsection{Unimodular random graphs}
\label{unimod}

The Mass-Transport Principle has been introduced by H\"aggstr\"om in \cite{Hag97} to study percolation and was further developed in \cite{BLPS99}. A random rooted graph $(G,\rho)$ obeys the \emph{Mass-Transport principle}  (abbreviated by MTP)  if for every Borel positive function $F : \mathcal{G}_{\bullet\bullet}  \to \R_{+}$ we have
\begin{eqnarray} \mathbf{E}\left[\sum_{x \in \op{V}(G)} F(G,\rho,x) \right ] &=&  \mathbf{E}\left[\sum_{x\in \op{V}(G)} F(G,x,\rho)\right] \label{MTP}. \end{eqnarray}
The name comes from the interpretation of $F$ as an amount of mass sent from $\rho$ to $x$ in $G$: the mean amount of mass that $\rho$ receives is equal to the mean quantity it sends. The MTP holds for a great variety of random graphs,  see \cite{AL07} where the MTP is extensively studied.

\begin{definition}{\cite[Definition 2.1]{AL07}} If $(G,\rho)$ satisfies \eqref{MTP} it is called \emph{unimodular} (See \cite{AL07} for explanation of the terminology).
\end{definition}

Let us explain the link between unimodular random graphs  and reversible random graphs.
Suppose that $F : \mathcal{G}_{\bullet\bullet} \to \R_{+}$ is a Borel positive function such that $F$ is supported by the subset of $ \mathcal{G}_{ \bullet\bullet}$ determined by the condition that the roots are neighbors, that is
\begin{eqnarray}
\label{form} F(G,x,y) &=&F(G,x,y)\mathbf{1}_{x\sim y}. \end{eqnarray}
Applying the MTP to a unimodular random graph $(G,\rho)$ with the function $F$ we get \begin{eqnarray} \mathbf{E}\left[\sum_{x\sim \rho} F(G,\rho,x)\right] &=& \mathbf{E}\left[\sum_{x\sim \rho} F(G,x,\rho)\right],
 \nonumber  \end{eqnarray}
 or equivalently
 \begin{eqnarray} \mathbf{E}\left[\op{deg}(\rho) \frac{1}{\op{deg}(\rho)}\sum_{x\sim \rho} F(G,\rho,x)\right] &=& \mathbf{E}\left[\op{deg}(\rho) \frac{1}{\op{deg}(\rho)}\sum_{x\sim \rho} F(G,x,\rho)\right].
 \nonumber  \end{eqnarray}

Let $(\tilde{G},\tilde{\rho})$ be distributed according to $(G,\rho)$ biased by $\op{deg}(\rho)$ (assuming that $\mathbf{E}\left[\op{deg}(\rho)\right] < \infty$), that is for any Borel $f : \mathcal{G}_{\bullet} \to \mathbb{R}_{+}$ we have $ \mathbf{E}[f(\tilde{G},\tilde{\rho})]= \mathbf{E}[ \mathrm{deg}(\rho)]^{-1} \mathbf{E}[f(G,\rho)  \mathrm{deg}(\rho)]$.  Conditionally on $(\tilde{G},\tilde{\rho})$, let $X_{1}$ be a one-step simple random walk starting on $\tilde{\rho}$ in $\tilde{G}$. Then the last display  is equivalent to 
\begin{eqnarray} \label{starev} (\tilde{G},\tilde{\rho},X_{1}) &\overset{(d)}{=}& (\tilde{G},X_{1},\tilde{\rho}).\end{eqnarray}
The graph $(\tilde{G},\tilde{\rho})$ is thus reversible hence stationary.
Reciprocally, if $(\tilde{G},\tilde{\rho})$ is reversible we deduce that the graph $(G,\rho)$ obtained after biasing by $\op{deg}(\rho)^{-1}$ obeys the MTP with functions of the form $F(G,x,y) \mathbf{1}_{x\sim y}$.  By \cite[Proposition 2.2]{AL07} this is sufficient to imply the full mass transport principle. Let us sum-up.
\begin{proposition} \label{unista} There is a correspondence between unimodular random graphs such that the expectation of the degree of the root is finite and reversible random graphs:
\begin{eqnarray*}
 (G,\rho) \mbox{ unimodular and } \mathbf{E}[\op{deg}(\rho)] < \infty &\overset{\mbox{\footnotesize bias by } \op{deg}(\rho)}{\underset{\mbox{\footnotesize bias by }\op{deg}(\rho)^{-1}}{\rightleftarrows}}& (G,\rho) \mbox{ reversible}.
 \end{eqnarray*}
 \end{proposition}

\subsection{Measured equivalence relations}
In this section we recall the notion of measured graphed equivalence relation. This concept will not   be used in the rest of the paper.\\

Let $(B,\mu)$ be a standard Borel space with a probability measure $\mu$ and let ${E} \subset B^2$ be a symmetric Borel set. We denote the smallest equivalence relation containing $E$ by $\mathcal{R}$. Under mild assumptions the triplet $(B,\mu,E)$ is called \emph{a measured graphed equivalence relation} (MGER). The set ${E}$ induces a graph structure on $B$ by setting $x\sim y \in B$ if $(x,y)\in {E}$. For $x \in B$, one can interpret the equivalence class of $x$ as a graph with the edge set given by $E$, which we root at the point $x$. If $x$ is sampled according to $\mu$, any measured graphed equivalence relation can be seen as a random rooted graph. See \cite{AL07,Kai98,Kai01,KS10,Pau99}. 

Reciprocally, the Polish space $\mathcal{G}_{\bullet}$ can be equipped with a symmetric Borel set $E$ where $((G,\rho), (G',\rho')) \in {E}$ if $(G,\rho)$ and $(G',\rho')$ represent the same isomorphism class of non-rooted graphs but are rooted at two different neighboring vertices. Denote $\mathcal{R}$ the smallest equivalence relation on $\mathcal{G}_{\bullet}$ that contains $E$. Thus a random rooted graph $(G,\rho)$ of distribution $\mathbf{P}$ gives rise to $(\mathcal{G}_{\bullet},\mathbf{P},E)$ which, under mild assumptions on $(G,\rho)$ is a measured graphed equivalence relation. \\

Remark however that the measured graphed equivalence relation we obtain with this procedure can have a graph structure on equivalence classes very different from what we could expect : Consider for example the (random) graph $\Z^2$ rooted at $(0,0)$. Since $\Z^2$ is a transitive graph, the measure obtained on $\mathcal{G}_{\bullet}$ by the above procedure is concentrated on the singleton corresponding to the isomorphism class of $(\Z^2,(0,0))$. Hence the random graph associated to this MGER is the rooted graph with one point, which is quite different from the graph $\Z^2$ we could expect!

 There are two ways to bypass this difficulty: considering \emph{rigid} graphs (that are graphs without non trivial isomorphisms see \cite[Section 1E]{KS10}) or add independent uniform labels $\in [0,1]$ on the graphs (see \cite[Example 9.9]{AL07}). Both procedures yield a MGER whose graph structure is that of $(G,\rho)$.\\

In particular we have the following correspondence between the notions of harmonic MGER and stationary random graph, totally invariant MGER  and reversible random graph, measure preserving MGER or and unimodular random graph, see \cite{AL07,Gab05,Kai98,Kai01,Pau99}. Also, the entropy theory has been developed in the context of random walks on equivalence relations, see \cite{CCC11,Bow10,Kai98,Kai01}. In the next section we will develop it, from scractch,  in the context of stationary random graphs.

\section{The Liouville property}In this section, we extend a well-known result on groups first proved in \cite{Ave74} relating Poisson boundary to entropy of a group. Here we adapt the proof which was given in \cite[Theorem 1]{KV83} in the case of groups (see also \cite{KW02} in the case of homogeneous graphs).  We basically follow the argument of \cite{KV83} using \emph{expectation} of entropy. The stationarity of the underlying random graph together with the Markov property of the simple random walk will replace homogeneity of the graph. We introduce the \emph{mean entropy} of the random walk and prove some useful lemmas. Then we derive the main results of this section. \medskip

In the following $(G,\rho)$ is a stationary random graph. Recall that conditionally on $(G,\rho)$, $\mathrm{P}_{x}^G$ is the law of the simple random walk $(X_n)_{n \geq 0}$ on $G$  starting from $x\in \op{V}(G)$. For every integer $0 \leq a \leq b < +\infty$, the \emph{entropy} of the simple random walk started at $x\in \op{V}(G)$ between times $a$ and $b$ is
$$ H_{a}^b(G,x) =  \sum_{x_{a}, x_{a+1}, \ldots , x_{b}} \varphi\left(\mathrm{P}_{x}^G(X_{a}=x_{a}, \ldots ,X_{b}=x_{b})\right),$$
where $\varphi(t) = -t\log(t)$. To simplify notation we write $H_{a}(G,x) =H_{a}^a(G,x)$. Recalling that $(G,\rho)$ is a random graph we set
$$ h_{a}^b = \mathbf{E}\left[H_{a}^b(G,\rho)\right] \mbox{ and } h_{a} = \mathbf{E}\left[H_{a}(G,\rho)\right].$$

\begin{proposition} If $(G,\rho)$ is a stationary random graph then $(h_{n})_{n\geq 0}$ is a subadditive sequence.
\end{proposition}
\proof Let $n,m \geq 0$. We have
\begin{eqnarray*} H_{n+m}(G,\rho) &=& \sum_{x_{n+m}} \varphi\left(\mathrm{P}_{\rho}^{G}(X_{n+m}=x_{n+m})\right).
\end{eqnarray*}
Applying the Markov property at time $n$, we get \begin{eqnarray*}
H_{n+m}(G,\rho) &=& \sum_{x_{n+m}}\varphi\left( \sum_{x_{n}}\mathrm{P}_{\rho}^G(X_{n}=x_{n})\mathrm{P}_{x_{n}}^G(X_{m}=x_{n+m})\right).\end{eqnarray*} Since $\varphi$ is concave and $\varphi(0)=0$  we have $\varphi(x+y) \leq \varphi(x)+\varphi(y)$,for every $x,y \geq 0$. Hence we obtain
\begin{eqnarray*} H_{n+m}(G,\rho)  &\leq & \sum_{x_{n+m}}\sum_{x_{n}}\varphi\left(\mathrm{P}_{\rho}^G(X_{n}=x_{n})\mathrm{P}_{x_{n}}^G(X_{m}=x_{n+m})\right)\\
  & = & H_{n}(G,\rho) + \sum_{x_{n}}\mathrm{P}_{\rho}^G(X_{n}=x_{n}) H_{m}(G,x_{n}).
 \end{eqnarray*}
  Taking expectations one has using (\ref{sta})
\begin{eqnarray} h_{n+m} &\leq& h_{n}+\mathbf{E}\left[\sum_{x_{n}}  \mathrm{P}_{\rho}^G(X_{n}=x_{n})H_{m}(G,x_{n})\right] \nonumber \\
 &=& h_{n} + \mathbb{E}\left[H_{m}(G,X_{n})\right] \quad = \quad h_{n}+h_{m}. \nonumber \end{eqnarray}\endproof
The subadditive lemma then implies that
 \begin{eqnarray}
 \label{h} \frac{h_{n}}{n} \tends{}{n\to \infty}h \geq 0.\end{eqnarray}  This limit is called the \emph{mean entropy} of the stationary random graph $(G,\rho).$  It plays the role of the (deterministic) entropy of a random walk on a group. In the rest of the paper, we will assume that $h$ is finite. The following theorem generalizes the well-known connection between the Liouville property and the entropy.

  \begin{theorem} \label{Liouville}  Let $(G,\rho)$ be a stationary random graph. The following conditions are equivalent:\begin{itemize} \item the tail $\sigma$-algebra associated to the simple random walk on $G$ started from $\rho$ is almost surely trivial (in particular it implies that $(G,\rho)$ is almost surely Liouville),
\item the mean entropy $h$ of $(G,\rho)$ is null.
\end{itemize}
\end{theorem}
Before doing the proof, we start with a few lemmas.

 \begin{lemma} \label{lemma2}For every $0\leq a \leq b < \infty$ we have
 $ h_{a}^b = h_{a} + (b-a)h_{1}.$ In particular for $k\geq 1$ we have $h_{1}^k=kh_{1}$.
 \end{lemma}
 \proof Let $0 \leq a \leq b <\infty$. An application of the Markov property at time $a$ leads to
 \begin{eqnarray*}
 H_{a}^b(G,\rho) &=& - \sum_{x_{a}, \ldots , x_{b}}\mathrm{P}_{\rho}^G(X_{a}=x_{a}, \ldots , X_{b}=x_{b})\log\left(\mathrm{P}_{\rho}^G(X_{a}=x_{a}, \ldots , X_{b}=x_{b}) \right)\\
 & = & - \sum_{x_{a}} \mathrm{P}_{\rho}^G(X_{a}=x_{a})\log\left(\mathrm{P}_{\rho}^G(X_{a}=x_{a})\right) + \sum_{x_{a}} \mathrm{P}_{\rho}^G(X_{a}=x_{a}) H_{1}^{b-a}(G,x_{a}). \end{eqnarray*}
  Taking expectations we get $h_{a}^b = h_{a} + h_{1}^{b-a}$. An iteration of the argument proves the lemma.
 \endproof
If $(G,\rho)$ is a fixed rooted graph and $(X_{n})_{n\geq0}$ is distributed according to $\mathrm{P}_{\rho}^G$, we introduce the following $\sigma$-algebra:
 \begin{eqnarray*}
 \mathcal{F}_{n}(G,\rho) &=& \sigma(X_{1}, \ldots ,X_{n}),\\
 \mathcal{F}^n(G,\rho) &=& \sigma(X_{n}, \ldots ),\\
 \mathcal{F}^\infty(G,\rho) &=& \bigcap_{n\geq 0} \mathcal{F}^n(G,\rho).
 \end{eqnarray*}

 The elements of the last $\sigma$-algebra are called tail events. By classical results of entropy theory, for all $k \geq 0$, the conditional entropy $H(\mathcal{F}_{k}(G,\rho)\mid\mathcal{F}^n(G,\rho))$ increases as $n\to \infty$ and converges to $H(\mathcal{F}_{k}(G,\rho)\mid\mathcal{F}^{\infty}(G,\rho))$. Furthermore, we have
 $$ H(\mathcal{F}_{k}(G,\rho)\mid\mathcal{F}^{\infty}(G,\rho)) \leq H(\mathcal{F}_{k}(G,\rho)),$$
 with equality if and only if $\mathcal{F}_{k}(G,\rho) \mbox{ and } \mathcal{F}^{\infty}(G,\rho)$ are independent. 

 \begin{lemma} For $1\leq k \leq n \leq m < +\infty$ we have
 $ \mathbf{E}\left[H(X_{1},\ldots, X_{k}\mid X_{n}, \ldots , X_{m})\right]= kh_{1}+h_{n-k}-h_{n}.$
 \end{lemma}
 \proof We have by definition
 \begin{eqnarray*}
& & H(X_{1}, \ldots, X_{k}\mid X_{n}, \ldots , X_{m})\\
 &= & -\sum_{\begin{subarray}{l}x_{1}, \ldots , x_{k} \\ x_{n}, \ldots , x_{m} \end{subarray}} \mathrm{P}_{\rho}^G(X_{i}=x_{i}, 1\leq i \leq k \mbox{ and } n\leq i \leq m) \log \left( \frac{\mathrm{P}_{\rho}^G(X_{i}=x_{i}, 1\leq i \leq k \mbox{ and } n\leq i \leq m)}{\mathrm{P}_{\rho}^G(X_{i}=x_{i}, n\leq i \leq m)}\right).
 \end{eqnarray*}
 Applying the Markov property at time $k$ one gets
 \begin{eqnarray*}
 & = & H_{1}^k(G,\rho)-H_{n}^m(G,\rho) + \sum_{x_{k}}\mathrm{P}_{\rho}^G(X_{k}=x_{k}) H_{n-k}^{m-k}(G,x_{k}),
 \end{eqnarray*}
 and taking expectations using \eqref{sta}, the right-hand side becomes $h_{1}^k-h_{n}^m+h_{n-k}^{m-k}$. An application of Lemma \ref{lemma2} completes the proof.
 \endproof
 \noindent In particular we see that the expected value of $H(X_{1}, \ldots , X_{k}\mid X_{n}, \ldots , X_{m})$  does not depend upon $m$ (this is also true without taking expectation and follows from Markov property at time $n$). If we let $m \to \infty$ in the statement of the last lemma, we get by monotonicity of conditional entropy and monotone convergence
\begin{equation} \label{main} \mathbf{E}\left[H(\mathcal{F}_{k}(G,\rho)\mid\mathcal{F}^n(G,\rho))\right] = kh_{1}+h_{n-k}-h_{n}.\end{equation}
\proof[Proof of Theorem \ref{Liouville}] Using again the monotonicity of conditional entropy $$H(\mathcal{F}_{1}(G,\rho) \mid \mathcal{F}^n(G,\rho)) \leq H(\mathcal{F}_{1}(G,\rho) \mid \mathcal{F}^{n+1}(G,\rho))$$
and the equality \eqref{main} for $k=1$, we deduce that
$(h_{n+1}-h_{n})_{n\geq 0}$ is decreasing  and converges towards $\tilde{h} \geq 0$. By (\ref{h}) and Cesaro's Theorem, we deduce that $\tilde{h}=h$. Thus letting $n \to\infty$ in (\ref{main}) we get by monotone convergence
$$ \mathbf{E}\left[H(\mathcal{F}_{k}(G,\rho)\mid\mathcal{F}^\infty(G,\rho))\right] = k(h_{1}-h).$$
 Comparing the last display with Lemma \ref{lemma2} (note that  $H(\mathcal{F}_{k}(G,\rho))= H_{1}^k(G,\rho)$), it follows that $h=0$ if and only if almost surely, for all $k \geq 0$, $\mathcal{F}^\infty(G,\rho)$ is independent of $\mathcal{F}_{k}(G,\rho)$.  Since there are no non-trivial events independent of all the coordinate $\sigma$-algebra we deduce that $\mathcal{F}^\infty(G,\rho)$ is almost surely trivial, in particular $(G,\rho)$ is Liouville. This completes the proof of Theorem \ref{Liouville}. \endproof

\proof[Proof of Theorem \ref{subexp}] Let $(G,\rho)$ be a stationary random graph of \emph{subexponential growth} that is
 $\mathbf{E}[\log(\#B_{G}(\rho,n))] = o(n),$  as $n \to \infty$.
Thanks to Theorem \ref{Liouville}, we only have to prove that the mean entropy of $G$ is zero. But by a classical inequality we have
$  H_{n}(G,\rho) \leq \log( \# B_{G}(\rho,n)),$ taking expectations and using \eqref{h} yields the result.  \endproof

In the preceding theorem we saw that subexponential growth plays a crucial role. In the case of transitive or Cayley graphs, all the graphs considered have at most an exponential growth. But  there are stationary graphs with superexponential growth, here is an example.
\begin{exe} Let $(G,\rho)$ be an augmented Galton-Watson tree (see Example \ref{AGW}) with offspring distribution  $(p_{k})_{k\geq 1}$ such that $\sum_{k\geq 1}k p_{k}= \infty$.  We have $$
\liminf_{n \to \infty} \frac{\mathbf{E}[\log\left(B_{G}(\rho,n)\right)]}{n} = \infty.$$
 \end{exe}

We can also extend the ``fundamental inequality'' for groups \cite{Gui80} or homogeneous graphs \cite{KW02}. The proof is mutatis-mutandis the same as in the group case

\begin{proposition} \label{fundamental} Let $(G,\rho)$ be a stationary and ergodic random graph of degree almost surely bounded by $M>0$. Conditionally on $(G,\rho)$, let $(X_{n})_{n\geq 0}$ be the simple random walk on $G$ starting from $\rho$. We denote \emph{the speed}  of the random walk by $s$ and \emph{the exponential volume growth}  of $G$ by $v$, namely \begin{eqnarray*} s &=& \limsup_{n\to \infty} n^{-1}\E{\op{d}_{\op{gr}}^G(X_{0},X_{n})},\\ v&=& \limsup_{n\to\infty} n^{-1}\Eb{\log(\# B_{G}(\rho,n))}.\end{eqnarray*} Then the mean entropy $h$ of $(G,\rho)$ satisfies $$
\frac{s^2}{2} \ \leq \   h \ \leq \ vs.$$ In particular $h =0 \iff s =0$ and if $s$ or $v$ is null then $(G,\rho)$ is almost surely Liouville. \end{proposition} 
\proof  Since $(G,\rho)$ is ergodic, we know from Theorem \ref{range}\eqref{speed} that $n^{-1} \op{d}_{\op{gr}}^G(X_0,X_n)$ converges almost surely and in $\mathbb{L}^1(\mathbb{P})$ towards $s \geq0$. In particular if $s>0$, for every $\varepsilon\in ]0,s[$ we have \begin{eqnarray} \mathbb{P}\big((s-\varepsilon)n \leq\op{d}_{\op{gr}}^G(X_{0},X_{n}) \leq (s+\varepsilon) n \big) &\underset{n \to \infty}{\longrightarrow} & 1.\label{s-+} \end{eqnarray}
\noindent\textit{Lower bound.} We suppose $s>0$ otherwise the lower bound is trivial. We have \begin{eqnarray}
H_{n}(G,\rho) &\geq& \sum_{\begin{subarray}{c} x_{n} \\ \op{d}_{\op{gr}}^G(\rho,x_{n}) \geq (s-\varepsilon)n \end{subarray}} \varphi(\mathrm{P}_{\rho}^G (X_{n}=x_{n})) \nonumber \\
& = & -\sum_{\begin{subarray}{c} x_{n} \\ \op{d}_{\op{gr}}^G(\rho,x_{n}) \geq (s-\varepsilon)n \end{subarray}} \mathrm{P}_{\rho}^G(X_{n}=x_{n}) \log\left(\mathrm{P}_{\rho}^G(X_{n}=x_{n}) \right) \nonumber \end{eqnarray}
At this point we use the Varopoulos-Carne estimates (see \cite[Theorem 12.1]{LP10}), for the probability inside the logarithm. Hence,
 \begin{eqnarray} H_{n}(G,\rho)& \geq & -\sum_{\begin{subarray}{c} x_{n} \\ \op{d}_{\op{gr}}^G(\rho,x_{n}) \geq (s-\varepsilon)n \end{subarray}} \mathrm{P}_{\rho}^G(X_{n}=x_{n}) \log\left( 2\sqrt{M} \exp\left( - \frac{(s-\varepsilon)^2n}{2}\right) \right) \nonumber \\
 & = & \log\left( 2\sqrt{M} \exp\left( - \frac{(s-\varepsilon)^2n}{2}\right) \right) \mathrm{P}_{\rho}^G\big(\op{d}_{\op{gr}}^G(X_{0},X_{n}) \geq (s-\varepsilon)n\big).\end{eqnarray}
 Now, we take expectation with respect to $\mathbf{E}$, divide by $n$ and let $n\to \infty$. Using \eqref{s-+} and \eqref{h} we have $h \geq \frac{(s-\varepsilon)^2}{2}$.

\noindent \textit{Upper bound.} Fix $\varepsilon >0$.  To simplify notation, we write $B_{s}$ for $B_{G}(\rho, (s+\varepsilon) n)$ and $B_{s}^c$ for $B_{G}(\rho,n) \backslash B_{G}(\rho,(s+\varepsilon)n)$. We decompose the entropy $H_{n}(G,\rho)$ as follows \begin{eqnarray*}  H_{n}(G,\rho) &=& \sum_{x_{n} \in B_{s}} \varphi(\mathrm{P}_{\rho}^G(X_{n}=x_{n})) + \sum_{x_{n} \in B_{s}^c } \varphi(\mathrm{P}_{\rho}^G(X_{n}=x_{n}))\\   & \leq & \left( \sum_{x_{n} \in B_{s}} \mathrm{P}_{\rho}^G (X_{n}=x_{n})\right)\log \left(  \frac{ \# B_{s}}{ \sum_{x_{n} \in B_{s}} \mathrm{P}_{\rho}^G(X_{n}=x_{n})} \right) \\   &+& \left(\sum_{x_{n} \in B_{s}^c } \mathrm{P}_{\rho}^G(X_{n}=x_{n}) \right) \log\left( \frac{\# (B_{s}^c )}{\sum_{x_{n} \in B_{s}^c } \mathrm{P}_{\rho}^G(X_{n}=x_{n}) } \right). \\  \end{eqnarray*}  We used the concavity of $\varphi$ for the inequalities on the sums of the right-hand side.  Using the uniform bound on the degree, we get the crude upper bound $\# (B_{s}^c) \leq \# B_{G}(\rho,n) \leq M^n$. Taking expectation we obtain (using the easy fact that for $x \in [0,1]$ one has $-x\log(x) \leq e^{-1}$) $$ h_{n} \leq 2e^{-1}+ \Eb{\log\left(\# B_{G}(\rho,(s+\varepsilon)n)\right)} + \mathbb{P} \left( \op{d}_{\op{gr}}^G(X_{0},X_{n}) \geq (s+\varepsilon)n \right)n \log(M).$$ Divide the last quantities by $n$ and let $n \to \infty$, then \eqref{h} and \eqref{s-+} show that $h \leq (s+\varepsilon)v$. \endproof

\begin{rek}	 A natural question (raised by the referee) in this setting is whether $h=0$ is actually equivalent to the Liouville property. Also it would be nice to have a Shannon type convergence for the mean entropy, see \cite{KV83}. We did not pursue these goals herein. 
\end{rek}

 \section{The Radon-Nikodym Cocycle}
 In this part we borrow and reinterpret in probabilistic terms a notion coming from the measured equivalence relation theory, the Radon-Nikodym cocycle (see \cite{FM77}), in order  to deduce several properties of stationary non reversible graphs, (see e.g.\,\cite{KS10} for another application). This notion plays the role of the modular function in transitive graphs, see \cite{SW90}.  The results of this section are very close to known results in measured equivalence relation theory (see \cite{Kai98,Pau99}). Our emphasis being on the probabilistic interpretation of the Radon-Nikodym cocycle rather than on the results themselves. \medskip
 
   In the remainder of this section, $(G,\rho)$ is a stationary random graph whose degree is almost surely bounded by a constant $M>0$.\medskip 

 Conditionally on $(G,\rho)$ of law $\mathbf{P}$, let $(X_{n})_{n\geq0}$ be a simple random walk of law $\mathrm{P}_{\rho}^G$. Let $\mu_{\to}$ and $\mu_{\leftarrow}$ be the two probability measures on $\mathcal{G}_{\bullet\bullet}$ such that $\mu_{\to}$ is the law of $(G,X_{0},X_{1})$  and $\mu_{\leftarrow}$ that of $(G,X_{1},X_{0})$. It is easy to see that the two probability measures $\mu_{\to}$ and $\mu_{\leftarrow}$ are mutually absolutely continuous. To be precise, for any Borel set $A\subset\mathcal{G}_{\bullet\bullet}$, since $(G,\rho)$ is a stationary random graph \eqref{sta} we have
  \begin{eqnarray*}\mathbb{P}\left((G,X_{0},X_{1}) \in A\right) &=&  \mathbb{P}\left((G,X_{1},X_{2}) \in A\right) \\ & \geq & \mathbb{P}\left((G,X_{1},X_{0}) \in A\,,\ X_{2}=X_{0}\right)\\ 
  &\geq & M^{-1} \mathbb{P}\left((G,X_{1},X_{0}) \in A\right).  \end{eqnarray*}
 Thus the Radon-Nikodym derivative of $(G,X_{1},X_{0})$ with respect to $(G,X_{0},X_{1})$, given for any $(g,x,y) \in \mathcal{G}_{\bullet\bullet}$ such that $x \sim y$ by
 $$ \Delta(g,x,y) := \frac{\mathrm{d}\mu_{\leftarrow}}{\mathrm{d}\mu_{\to}}(g,x,y),$$ can be chosen such that \begin{eqnarray}
 M^{-1} \leq \Delta(g,x,y) \leq M. \label{bdd} \end{eqnarray}
   Note that the function $\Delta$ is defined up to a set of $\mu_{\to}$-measure zero, and in the following we fix an arbitrary representative satisfying \eqref{bdd} and we keep the notation $\Delta$ for this function. Since $\Delta$ is a Radon-Nikodym derivative we obviously have $\mathbb{E}[\Delta(G,X_{0},X_{1})]=1$ and  Jensen's inequality yields
\begin{eqnarray}
\mathbb{E}\Big[\log\big(\Delta(G,X_{0},X_{1})\big)\Big] &\leq &0, \label{jensenss}
\end{eqnarray}
with equality if and only if $\Delta(G,X_{0},X_{1})=1$ almost surely. In this latter case the two random variables $(G,X_{0},X_{1})$ and $(G,X_{1},X_{0})$ have the same law, that is $(G,\rho)$ is reversible.

   \begin{lemma}\label{countable} With the above notation, let $A$ be a Borel subset of $\mathcal{G}_{\bullet\bullet}$ of $\mu_{\to}$-measure zero. Then for $\mathbf{P}$-almost every rooted graph $(g,\rho)$ and every $x,y\in \op{V}(g)$ such that $x\sim y$ we have $(g,x,y) \notin A$.
   \end{lemma}
   \proof By stationarity, for any $n\geq 0$ the variable $(G,X_{n},X_{n+1})$ has the same distribution as $(G,X_{0},X_{1})$. Thus we have
   \begin{eqnarray} 0 &=& \sum_{n\geq 0} \mathbb{P}\left((G,X_{n},X_{n+1})\in A\right) = \E{\sum_{n\geq 0}\mathbf{1}_{(G,X_{n},X_{n+1})\in A}}\nonumber \\
   &=& \Eb{\sum_{x\sim y \in G}  \mathbf{1}_{(G,x,y) \in A}\left(\sum_{n\geq 0} \mathrm{P}_{\rho}^G(X_{n}=x, X_{n+1}=y)\right)}.\nonumber \end{eqnarray}

   Let $x\sim y$ in $G$. Since $G$ is connected, there exists values of  $n$ such that the probability that $X_{n}=x$ and $X_{n+1}=y$ is positive. Thus the sum between parentheses in the last display is positive. This proves the lemma.\endproof

 Note that the function $(g,x,y) \to \Delta(g,y,x)$ is also a version of the Radon-Nikodym derivative $\frac{\mathrm{d}\mu_{\to}}{\mathrm{d}\mu_{\leftarrow}}$, hence we have $\Delta(g,x,y) = \Delta(g,y,x) ^{-1}$ for $\mu_{\to}$-almost every bi-rooted graphs in $\mathcal{G}_{\bullet\bullet}$. By the above lemma we also have $\Delta(g,x,y)= \Delta(g,y,x)^{-1}$ for $\mathbf{P}$-almost every rooted graph $(g,\rho)$ and every vertices $x,y \in \op{V}(g)$ such that $x\sim y$.

\begin{lemma} \label{welldefined} For $\mathbf{P}$-almost every $(g,\rho)$, and every cycle $\rho = x_{0} \sim x_{1} \sim \ldots \sim x_{n} = \rho$ in $g$ we have
\begin{eqnarray} \label{cycle} \prod_{i=0}^{n-1} \Delta(g,x_{i},x_{i+1}) &=&1. \end{eqnarray}
\end{lemma}
 \proof 
By a standard calculation on the simple random walk, conditionally on $(G,\rho)$ and on $\{ \rho=X_{0}= X_{n}\}$, the path $(X_{0}, X_{1}, \ldots ,X_{n-1}, X_{n})$ has the same distribution as the reversed one $(X_{n}, X_{n-1}, \ldots ,X_{1}, X_{0})$. In other words, for any positive Borel function $F : \R_{+} \to \R_{+}$ we have
 \begin{eqnarray}
 \mathbb{E}\left[ F \left( \prod_{i=0}^{n-1} \Delta(G,X_{i},X_{i+1})\right) \mathbf{1}_{X_{n}=X_{0}}\right] & =&   \mathbb{E}\left[ F \left( \prod_{i=0}^{n-1} \Delta(G,X_{i+1},X_{i})\right)\mathbf{1}_{X_{n}=X_{0}}\right]\nonumber \\
 &=&  \mathbb{E}\left[ F \left( \prod_{i=0}^{n-1} \Delta(G,X_{i},X_{i+1})^{-1}\right)\mathbf{1}_{X_{n}=X_{0}}\right]. \nonumber\end{eqnarray}
 Where we used the fact that for $\mathbf{P}$-almost every $(g,\rho)$ and  for any neighboring vertices $x,y \in \op{V}(g)$, we have $\Delta(g,x,y) = \Delta(g,y,x)^{-1}$.  Since for every $(g,\rho)\in \mathcal{G}_{\bullet}$ and any cycle $\rho=x_{0}\sim x_{1} \sim \ldots \sim x_{n}=\rho$ we have $P_{\rho}^G(X_{0}=x_{0}, X_{1}=x_{1}, \ldots , X_{n}=x_{n})>0$ the desired result easily follows. \endproof

 Suppose that the above lemma hold, then we can extend the definition of $\Delta$ to an arbitrary (isomorphism class of) bi-rooted graph $(g,x,y)$ without assuming $x\sim y$ (compared with \cite[Proof of Théorème 1.15 ]{Pau99}). If $x,y \in g$, let $x=x_{0} \sim x_{1} \sim \ldots \sim x_{n} =y$ be a path in $g$ between $x$ and $y$, and set

\begin{eqnarray}  \Delta(g,x,y) &:=& \prod_{i=0}^{n-1} \Delta(g,x_{i},x_{i+1}), \label{defrn}
 \end{eqnarray}
 and by convention $\Delta(g,x,x)=1$.
 This definition does not depend on the path chosen from $x$ to $y$ by the last lemma and is well founded for $\mathbf{P}$-almost every graph $(g,\rho)$ and every $x,y \in \op{V}(g)$.  We can now prove (compare with \cite{Pau99}):
 \begin{theorem}Let $(G,\rho)$ be a stationary ergodic random graph. Assume that $(G,\rho)$ is not reversible. Then almost surely the function $$x \in \op{V}(G) \mapsto \Delta(G,\rho, x),$$ is positive harmonic and non constant. In particular $(G,\rho)$ is almost surely transient.
 \end{theorem}
 \proof We follow the proof of \cite{Pau99}. By the stationarity of $(G,\rho)$, for any Borel function $F:\mathcal{G}_{\bullet} \to \R_{+}$ we have
 \begin{eqnarray} \E{F(G,X_{0})}  =  \E{F(G,X_{1})} = \E{F(G,X_{0}) \Delta(G,X_{0},X_{1})}. \nonumber \end{eqnarray}
 We thus get 
$\op{deg}(\rho)^{-1}\sum_{\rho \sim x} \Delta(G,\rho,x)= 1$ almost surely. 
 It follows from Lemma \ref{countable}, that almost surely, for any $x \in \op{V}(G)$ we have $$\frac{1}{\op{deg}(x)}\sum_{x \sim y} \Delta(G,x,y) =1.
$$ One gets from the previous display and the definition of $\Delta$, that $x \mapsto \Delta(G,\rho, x)$ is almost surely harmonic. Notice that if $x \mapsto \Delta(G,\rho,x)$ is constant then this constant is $1$. Also, the event $\{ \Delta(G,.,.) \mbox{ is constant}\}$ is an event which is invariant by the shift under $ \mathbb{P}$, more precisely if $\Delta(G,.,.)$ is constant over $(G,\rho)$ it is also constant over $(G,X_1)$. Thus, by ergodicity if $x \mapsto \Delta(G, \rho,x)$ has a positive probability to be constant then it is almost surely constant, and this constant equals $1$. This case is excluded because $(G,\rho)$ is not reversible. By standard properties of random walk on graphs, the existence of a non constant positive harmonic function implies transience. \endproof

With the hypotheses of the preceding theorem, we have in fact much more than transience of almost every graph $(G,\rho)$: The simple random walk is ballistic! This phenomenon has been known for long in the context of foliations, see \cite{Kai88}.

 \begin{theorem} \label{ballistic}Let $(G,\rho)$ be a stationary and ergodic random graph of degree almost surely bounded by $M>0$. If $(G,\rho)$ is non reversible, then the speed $s$ (see \eqref{speed}) of the simple random walk on $(G,\rho)$  is positive.
 \end{theorem}
 \proof The idea is to consider the rate of growth of the Radon-Nikodym cocycle along sample paths as in \cite[Corollary 1 of Theorem 2.4.2]{Kai98}. We consider the random process $(\log(\Delta(G,X_{0},X_{n})))_{n\geq0}$. By Proposition \ref{welldefined} we almost surely have for all $n\geq 0$
 \begin{eqnarray}
 \log\big(\Delta(G,X_{0},X_{n})\big) & =& \sum_{i=0}^{n-1}  \log\big(\Delta(G,X_{i},X_{i+1})\big). \end{eqnarray}
By \eqref{bdd} we have $\mathbb{E}[|\log(\Delta(G,X_{0},X_{1}))|] < \infty$ and the ergodic theorem implies the following almost sure and $\mathbb{L}^1$ convergence with respect to $\mathbb{P}$
 \begin{eqnarray} \label{ergogo} \frac{\log\big(\Delta(G,X_{0},X_{n})\big)}{n} &\underset{n\to \infty}{\longrightarrow} &\mathbb{E}[\log(\Delta(G,X_{0},X_{1}))].\end{eqnarray}
By computing $\Delta(G,X_{0},X_{n})$ as in \eqref{defrn} along a geodesic path from $X_{0}$ to $X_{n}$ in $G$ and using \eqref{bdd} we deduce that a.s. for every $n\geq0$ $$|\log(\Delta(G,X_{0},X_{n}))| \leq \log(M)\op{d}_{\op{gr}}^G(X_{0},X_{n}).$$ If $(G,\rho)$ is not reversible, we already noticed that the inequality \eqref{jensenss} is strict. Thus combining \eqref{speed},\eqref{ergogo} and the last display we get  
$
s \geq|\mathbb{E}[\log(\Delta(G,X_{0},X_{1}))]|\log(M)^{-1}>0$, which is the desired result.\endproof

\begin{rek} By Corollary \ref{fundamental}, subexponential growth in the sense of \eqref{subexpo} implies $s=0$ for stationary and ergodic random graphs of bounded degree, so in particular such random graphs are reversible (known in the case of foliations, see \cite{Kai88}). This fact also holds without the bounded degree assumption (Russell Lyons, personal communication).

 \end{rek}
\section{Applications}

\subsection{The Uniform planar quadrangulation}

A planar map is an embedding of a planar graph into the two-dimensional sphere  seen up to continuous deformations. A quadrangulation is a planar map whose faces all have degree four. The Uniform Infinite Planar Quadrangulation (UIPQ)  introduced by Krikun in \cite{Kri05} is the weak local limit (in a sense related to $\op{d_{loc}}$) of uniform quadrangulations with $n$ faces with a distinguished oriented edge (see Angel and Schramm \cite{AS03} for previous work on triangulations).  We will not discuss the subtleties of planar maps nor the details of the construction of the UIPQ and refer the interested reader to \cite{Kri05,LGM10,Men08}. \\

The UIPQ is a random infinite graph $\mathrm{Q}_{\infty}$ (which is viewed as embedded in the plane) given with a distinguished oriented edge $\vec{e}$. We will forget the planar structure of the UIPQ and get a random rooted graph $({Q}_{\infty},\rho)$, which is rooted at the origin $\rho$ of $\vec{e}$. This graph is stationary and reversible, see \cite{CMMinfini}. One of the main open questions about this random infinite graph is its conformal type, namely is it (almost surely) recurrent or transient? It has been conjectured in \cite{AS03} (for the related Uniform Infinite Planar Triangulation) that ${Q}_{\infty}$ is almost surely recurrent. Although we know that the conformal type of the Riemann surface obtained from the UIPQ by gluing squares along edges  is  parabolic \cite{GR10} (see \cite{AS03} for related result on the Circle Packing), yet the absence of the bounded degree property prevents one from using the results of \cite{BS01} to get recurrence of the simple random walk on the UIPQ. Corollary \ref{uipq} may be seen as providing a first step towards the recurrence proof.

\proof[Proof of Corollary \ref{uipq}] The random rooted graph $({Q}_{\infty},\rho)$ is a stationary and reversible random graph. A proof of this fact can be found in \cite{CMMinfini}. By virtue of Theorem \ref{subexp}, we just have to show that $({Q}_{\infty},\rho)$ is of subexponential growth. To be completely accurate, we have to note that the graph $(Q_\infty,\rho)$ is not simple, that is contains loops and multiple edges. However, it is easy to check that  Theorem \ref{subexp} still holds in this more general setting. Thanks to \cite{Men08}, we know that the random infinite quadrangulation investigated in \cite{CD06} has the same distribution as the UIPQ. Hence, the volume estimate of \cite{CD06} can be translated into  \begin{eqnarray} \label{estimate}\Eb{\#B_{{Q}_{\infty}}(\rho,n)} &=& \Theta(n^4). \end{eqnarray} Hence Jensen's inequality proves that the UIPQ is of subexponential growth in the sense of $\eqref{subexpo}$ which finishes the proof of the corollary.\endproof
This corollary does not use the planar structure of UIPQ but only the invariance with respect to SRW and the subexponential growth. We believe that the result of Corollary \ref{uipq} also holds for the UIPT. A detailed proof could be given along the preceding lines but would require an extension of the estimates \eqref{estimate} (Angel \cite{Ang03} provides almost sure estimates that are closely related to \eqref{estimate} for the UIPT).

\subsection{Long range percolation clusters}
Consider the graph obtained  from $\mathbb{Z}^d$ by adding an edge between each pair of distinct vertices $x,y \in \mathbb{Z}^d$ with probability $p_{x,y}$ independently of the other pairs. Assume that
$$ p_{x,y} = \beta|x-y|^{-s},$$ for some $\beta>0$ and $s>0$.
This model is called \emph{long range percolation}. Berger \cite{Berg02} proved in dimensions $d=1$ or $d=2$ that if $d<s<2d$, then conditionally on $0$ being in an infinite cluster, this cluster is almost surely transient. In the same paper the following question \cite[(6.3)]{Berg02} is addressed:
\begin{question} Are there nontrivial harmonic functions on the infinite cluster of long range percolation with $d<s<2d$ ?
\end{question}
We answer negatively this question for \emph{bounded} harmonic functions.
 \proof First we remark that by a general result (see \cite[Example 9.4]{AL07}), conditionally on the event that $0$ belongs to an infinite cluster $\mathcal{C}_{\infty}$, the random rooted graph $(\mathcal{C}_{\infty},0)$ is a unimodular random graph. Furthermore, since $s>d$ the expected degree of $0$ is finite. Hence, by Proposition \ref{unista}, the random graph $(\tilde{\mathcal{C}}_{\infty},\tilde{0})$ obtained by biasing $(\mathcal{C}_{\infty},0)$ with the degree of $0$ is stationary.  By Theorem \ref{subexp} it suffices to show that the graph $\tilde{\mathcal{C}}_{\infty}$ is of subexponential growth in the sense of \eqref{subexpo}. For that purpose, we use the estimates given in \cite[Theorem 3.1]{Bis09}. For $x \in \mathcal{C}_{\infty}$, denote the graph distance from $0$ to $x$ in $\mathcal{C}_{\infty}$ by $ \op{d}_{\op{gr}}^{\mathcal{C}_{\infty}}(0,x)$. Then  for each $s' \in (d,s)$ there are constants $c_{1}, c_{2} \in (0,+\infty)$ such that, for $\delta' = 1/ \log_{2}(2d/s')$,
$$\mathbf{P}\left( \op{d}_{\op{gr}}^{\mathcal{C}_{\infty}}(0,x) \leq n \right) \leq c_{1} \left( \frac{e^{c_{2}n^{1/\delta'}}}{|x|}\right)^{s'}.$$
In particular, we deduce that \begin{eqnarray}
\Eb{\#B_{\mathcal{C}_{\infty}}(0,n)} &\leq& \kappa_{1} \exp\left( \kappa_{2} n^{1/\delta'}\right), \label{eq1}\end{eqnarray} where $\kappa_{1}$ and $\kappa_{2}$ are positive constants. Remark that $\delta' >1$. Thus we have, if $\op{deg}(0)$ denotes the degree of $0$ in $\mathcal{C}_{\infty}$,
\begin{eqnarray} \mathbf{E}\big[\log ( \# B_{\tilde{\mathcal{C}}_{\infty}}(\tilde{0},n))\big] &= &\frac{1}{\mathbf{E}[\op{deg}(0)]}\mathbf{E}\big[\op{deg}(0)\log ( \# B_{\mathcal{C}_{\infty}}({0},n))\big]  \nonumber \\
 & \leq & \frac{1}{\mathbf{E}[\op{deg}(0)]} \sqrt{\mathbf{E}[\deg(0)^2]\mathbf{E}[\log^2(\# B_{\mathcal{C}_{\infty}}(0,n))]}, \label{eq2} \end{eqnarray}
by the Cauchy-Schwarz inequality.  Since $s>d$ it is easy to check that the second moment of $\op{deg}(0)$ is finite. Furthermore, the function $ x \mapsto \log^2(x)$ is concave on $]e,\infty[$ so by Jensen's inequality we have $$\Eb{\log^2(\#B_{\mathcal{C}_{\infty}}(0,n))} \leq \log^2\left(\Eb{\#B_{\mathcal{C}_{\infty}}(0,n)}+2\right).$$  Hence, combining the last display with \eqref{eq1} and \eqref{eq2} we deduce that $({\tilde{\mathcal{C}}}_{\infty},0)$ is of subexponential growth in the sense of \eqref{subexpo}. \endproof

\begin{rek} It is also possible to derive this corollary from \cite[Theorem 4]{Kai90}, however we preferred to stick to the context of unimodular random graphs. \end{rek}

Note that by similar considerations, clusters of any invariant percolation on a group, in which the clusters have subexponential volume growth are Liouville, see \cite{BLPS99} for many examples. This holds in particular for Bernoulli percolation on Cayley graphs of subexponential growth, e.g.  on $\mathbb{Z}^d$.

\subsection{Planarity}

Simply connected planar Riemannian surfaces are conformally equivalent either to the Euclidean or to
the hyperbolic plane. Thus they are either recurrent for Brownian motion or admit non constant bounded harmonic functions. The same alternative holds for planar graphs of \emph{bounded degree}. They are either recurrent for the simple random walk or admit non constant bounded harmonic functions \cite{BS96a}. Combining Theorem \ref{subexp} with these results related to  planarity yields:
\begin{corollary} \label{UIPQc}Let $(G,\rho)$ be a stationary random graph with subexponential growth in the sense of \eqref{subexpo}. Suppose furthermore that almost surely $(G,\rho)$ is planar and has bounded degree. Then $(G,\rho)$ is almost surely recurrent.
\end{corollary}
\proof We already know by Theorem \ref{subexp} that $(G,\rho)$ is almost surely Liouville. In \cite{BS96a} it is shown that a transient planar graph with bounded degree admits non constant bounded harmonic functions. Therefore $G$ must be recurrent almost surely.\endproof

Note that without the bounded degree assumption it is easy to construct planar transient Liouville graphs, see \cite{BS96a}. However these graphs are not stationary. The following construction shows that the bounded degree assumption is needed in the last corollary: We construct a stationary and reversible random graph that is of subexponential growth but transient (see Proposition \ref{counterexample}). \bigskip \paragraph{The example.} We consider the sequence $\epsilon_1, ... , \epsilon_n, ... \in \{1,2\}$ defined recursively as follows. Start with $\epsilon_1=1$, if $\epsilon_1, ... , \epsilon_k$ are constructed we let $\xi_k = \prod_{i=1}^k \epsilon_k$, and set $\epsilon_{k+1}=1$ if $\xi_k > k^4$ and $\epsilon_{k+1}=2$ otherwise. Clearly there exist constants $0<c<C<\infty$ such that $c k^4 \leq \xi_k \leq C k^4$ for every $k\geq1$. We now consider the tree $T_n$ of height $n$, starting from an initial ancestor at height $0$ such that each vertex at height $0 \leq k \leq n-1$ has $\epsilon_{n-k}$ children. Hence the tree $T_n$ has only simple or binary branchings. The \emph{depth} $D(u)$ of a vertex $u$ in $T_n$ is $n$ minus its height. For example, the leaves of $T_n$ have depth $0$. The depth of an edge is the maximal depth of its ends.

\begin{figure}[!h]
 \begin{center}
 \includegraphics[width=10cm]{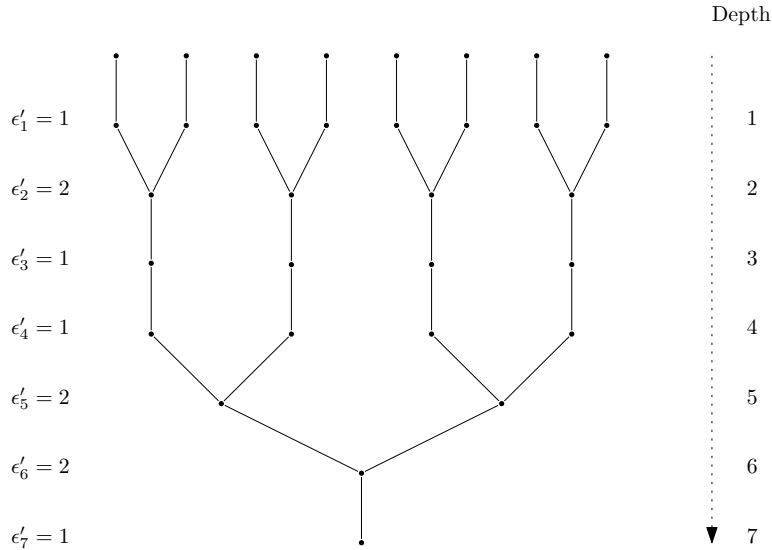}
 \caption{Construction of the tree $T'_7$ with a sample sequence $(\epsilon'_k)_{k\geq 1}$}
 \end{center}
 \end{figure}  We also introduce the infinite ``canopy" tree $T_\infty$ (the limit of the $T_n$'s seen from the top) which is a tree of infinite depth such that each vertex at depth $k$ is linked to $\epsilon_k$ vertices at depth $k-1$ for $k \geq 1$. \begin{lemma}  \label{growth}There exists a constant $C>0$ such that for every $u \in T_\infty$ and $r \geq 0$ we have 
  \begin{eqnarray*} \# B_{T_\infty}(u,r) & \leq & Cr^4.  \end{eqnarray*}
 \end{lemma}\proof Fix $r \in \{0,1,2, ...\}$. Let $u \in T_\infty$ and let $d\geq0$ be its depth. We suppose that $d > r$. We denote by $v$ the ancestor of $u$ at depth $r+d$. Clearly the ball of radius $r$ around $u$ in $T_\infty$ is contained in the subset of $T_\infty$ made of the vertices that are descendants of $v$ and whose depth is in between $d+r$ and $d-r$. This set has a cardinal equal to 
  \begin{eqnarray*} 1 + \epsilon_{d+r} + \epsilon_{d+r}\epsilon_{d+r-1} + ... + \epsilon_{d+r-1}\epsilon_{d+r-1}... \epsilon_{d-r-1} &=& \xi_{d+r}\left(\sum_{i=d-r}^{d+r} \xi_{i}^{-1}\right).  \end{eqnarray*}Recall that we have $ \xi_r = \Theta(r^4)$ as $r \to \infty$. Henceforth, when $d\leq 2r$, the last display is bounded from above by $\kappa (3r)^4 \sum_{1}^{\infty} \xi_k^{-1}$ and when $d \geq 2r$ we use the upper bound $\kappa' 2r ((r+d)/(r-d))^4$, where $\kappa,\kappa'>0$ are two positive constants independent of $d$ and $r$. In both cases we have $ \# B_{T_\infty}(u,r) = O(r^4)$. The case $r \geq d$ is similar and is left to the reader.\endproof 

Now we consider the graphs $T_n^{R}$ and $T_\infty^R$ obtained from $T_n$ and $T_\infty$ by replacing each edge at depth $k$ by $k^2$ parallel edges. The graph $T_\infty^R$ is obviously a tree with multiple edges that has only one end. We claim that this tree is transient, indeed its type is equivalent to that of a single spine with $k^2$ parallel edges at level $k$ (we can chop of the finite trees attached to the spine to study recurrence or transience). Since the conductance of the last spine is $ \sum_{k \geq 1} k^{-2} < \infty$ it is transient so is the tree $T_{\infty}^R$, see \cite{LP10}. 

\begin{figure}[!h]
 \begin{center}
 \includegraphics[width=13cm]{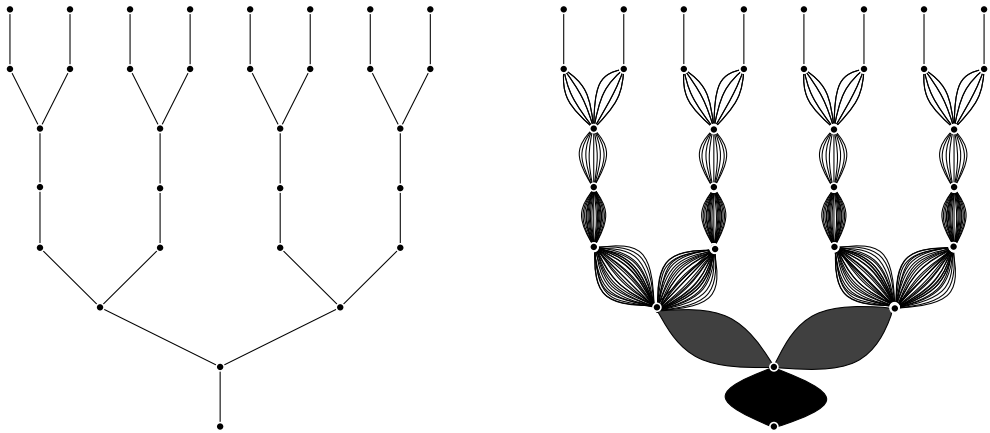}
 \caption{Muplication of edges.}
 \end{center}
 \end{figure}

We transform these deterministic graphs into random ones. The root $\rho_n$ is chosen among all vertices of $T_n^{R}$ proportionally to the degree. This boils down to picking an oriented edge uniformly at random in $T_n^R$ and consider its starting point $\rho_n$.
\begin{proposition} \label{counterexample} We have the convergence in distribution for $\op{d_{loc}}$
\begin{eqnarray} (T_n^R,\rho_n) &\underset{n \to \infty}{\longrightarrow} &(T^R_\infty,\rho),\end{eqnarray}
for a particular choice of a random root $\rho \in T_\infty^R$. In particular $(T^R_\infty,\rho)$ is a planar transient stationary and reversible random graph of subexponential growth.
\end{proposition}
\proof It is enough to show that $D(\rho_n)$ converges in distribution to a non degenerate random variable denoted by $D$ as $n\to \infty$. Indeed if we choose a random root $\rho \in T_\infty^R$ with depth given by $D$, since the $r$-neighborhood of a vertex at depth $k$ in $T_n^R$ and in $T_\infty^R$ are the same when $n \geq r+k$, we easily deduce the weak convergence of $(T_n^R,\rho_n)$ to $(T_\infty^R,\rho)$ for $ \op{d_{loc}}$. Furthermore since the random rooted graphs $(T_n^R,\rho_n)$ are stationary and reversible (see Example \ref{example}), the same holds for $(T_\infty^R, \rho)$ as weak limit  of stationary and reversible graphs in the sense of $\op{d_{loc}}$.\\
Let $k \geq 0$. The probability that $D(\rho_n)=k$ is exactly the proportion of oriented edges whose origin is a vertex of depth $k$. Thus with the convention $\xi_0,\xi_{-1}=1$ we have
$$ \mathbb{P}(D(\rho_n)=k) = \left(k^2 \frac{\xi_n}{\xi_{k-1}} + (k+1)^2 \frac{\xi_n}{\xi_{k}}\right)  \left(2\xi_n\sum_{i=0}^{n-1} \frac{(i+1)^2}{\xi_{i}} \right)^{-1}.$$Since $\xi_k \geq ck^4$, clearly the series $\sum i^{2}\xi_i^{-1}$ converges. Hence, the probabilities in the last display converge when $n\to \infty$, thus proving the convergence in distribution of $D(\rho_n)$. 
Furthermore, by Lemma \ref{growth}, $T_\infty^R$ is of subexponential growth. \endproof

\paragraph{Questions} \ \\

\noindent $\bullet$ In the preceding construction, the degree of $\rho$ in $T_{\infty}^R$ has a polynomial tail. Is it possible to construct a planar stationary and reversible graph of subexponential growth such that the degree of the root vertex has an exponential tail for which the SRW is transient?
\bigskip

\noindent $\bullet$ Let $(G,\rho)$ be a  limit in distribution of finite \emph{planar} stationary graphs for $\op{d_{loc}}$ (see \cite{BS01}). Is it the case that $(G,\rho)$ is almost surely Liouville\footnote{There are local limits of finite planar graphs with exponential growth. For example local limit of full binary trees up to level $n$ with the root picked according to the degree}? Does SRW on $(G,\rho)$ have zero speed? 

\bigskip

\noindent $\bullet$ In \cite{BCd-para} a generalization of limits of finite planar graphs to graphs associated to sphere packings in $\mathbb{R}^d$ was studied. Extend the preceding questions to these graphs. 

\bigskip

\noindent \textbf{Added in proof:} After the completion of this work, Gurel-Gurevich and Nachmias \cite{GGN12} proved that the UIPQ (and the UIPT) is recurrent which implies Corollary \ref{UIPQc}.\medskip 

\noindent \textbf{Acknowledgments :} We are grateful to Pierre Pansu, Frederic Paulin, Damien Gaboriau and Russell Lyons for many stimulating lessons on measured equivalence relations. We thank  Jean-Fran\c cois Le Gall for a careful reading of a first version of this paper. We are also indebted to Omer Angel for a discussion that led to Proposition \ref{counterexample}. Thanks also go to an anonymous referee for precious comments and references.

\end{document}